\documentclass{amsart}

\usepackage{amsmath}
\usepackage{amscd}
\usepackage{amssymb} 

\newcommand{\cal}{\mathcal}
\newcommand{\bk}{{\bf k}}

\newcommand{\bC}{{\Bbb C}}

\newcommand{\bQ}{{\Bbb Q}}
\newcommand{\bR}{{\Bbb R}}
\newcommand{\bZ}{{\Bbb Z}}
\newcommand{\cA}{{\cal A}}

\newcommand{\cD}{{\cal D}}

\newcommand{\cH}{{\cal H}}
\newcommand{\cL}{{\cal L}}

\newcommand{\fd}{{\frak d}}

\DeclareMathOperator{\ad}{ad}
\DeclareMathOperator{\End}{End}
\DeclareMathOperator{\Hom}{Hom}
\DeclareMathOperator{\Img}{Im}

\DeclareMathOperator{\Ker}{Ker}

\DeclareMathOperator{\vol}{vol}

\newtheorem{theorem}{Theorem}[section]
\newtheorem{theorem/definition}{Theorem/Definition}[section]

\newtheorem{lemma}{Lemma}[section]

\theoremstyle{remark}
\newtheorem{remark}{Remark}[section]

\theoremstyle{definition}
 \newtheorem{example}{Example}[section]

\begin{document}
\title
{Homological perturbation theory and mirror symmetry}
\author{Jian Zhou}
\address{Department of Mathematics\\
Texas A\&M University\\
College Station, TX 77843}
\email{zhou@math.tamu.edu}
\begin{abstract}
We explain how deformation theories of geometric objects such as
complex structures, Poisson structures and holomorphic bundle structures
lead to differential Gerstenhaber or Poisson algebras.
We use homological perturbation theory to obtain $A_{\infty}$ algebra
structures and some canonically defined deformations of such structures
on the cohomology.
We formulate the $A_{\infty}$ algebraic mirror symmetry as the identification
of the $A_{\infty}$ algebras together with their canonical deformations
constructed this way on different manifolds.
\end{abstract}
\maketitle

\section{Introduction}

We discuss in this paper
some applications of methods from algebraic topology to 
mirror symmetry. 
String theorists are interested in two kinds of conformal field
theories on any Calabi-Yau $3$-fold:
the $A$-type theories which do not depend on the complex structure on $M$
but is sensitive to the deformation of K\"{a}hler structure $M$,
the $B$-type theories which do not depend on the K\"{a}hler structure
but is sensitive to the deformation of the complex structure on $M$.
Roughly speaking,
the mirror symmetry conjecture states that 
for a suitable Calabi-Yau manifold $M$,
there exists another Calabi-Yau manifold $\widehat{M}$ such
that a specific $A$-type theory on $M$ can be identified with 
a specific $B$-type theory on $\widehat{M}$.
For more details,
see Yau \cite{Yau}. 
From the above stipulation on $A$-type theories and $B$-type theories,
it is conceivable that deformation theories of complex and 
K\"{a}hler (to be more precise, Poisson) structures lead
to a $B$-type theory and an $A$-type theory respectively.
In physics literature,
they correspond to the Kodaira-Spencer theory of gravity 
\cite{Ber-Cec-Oog-Vaf}
and the theory of K\"{a}her gravity respectively \cite{Ber-Sad}.

According to the well-known Deligne's principle:
``In characteristic zero a deformation problem
is controlled by a differential graded Lie algebra with the property that
quasi-isomorphic differential graded Lie algebras 
give the same deformation theory."
We emphasize that for many deformation problems
(including the aforementioned deformation problems 
on Calabi-Yau manifolds),
the controlling differential graded Lie algebra (DGLA)
is part of a differential Gerstenhaber, Gerstenhaber-Batalin-Vilkovisky
or Poisson algebra.
In other words,
by considering suitable extended deformation problems,
one encounters differential complexes 
with both a multiplicative structure and a Lie bracket 
which satisfy some compatibility conditions between them.
We will use homological perturbation theory to obtain $A_{\infty}$
structures  on the cohomology groups from the multiplicative structure,
and canonical deformations of these $A_{\infty}$ algebras 
using the Lie bracket.
We formulate mirror symmetry as the identifications of 
the resulting objects from different deformation problems.
Based on some earlier results in joint work with Cao 
\cite{Cao-Zho1, Cao-Zho2, Cao-Zho3},
and motivated by this work,
we gave in \cite{Cao-Zho4} a discussion of mirror symmetry
in term of formal Frobenius manifold structures 
by considering the DGBV algebras arising from 
the extended deformation problems  
of the complex structure and Poisson structure on a Calabi-Yau manifold.
We will compare that approach with the approach in this work
in a separate paper.

Our formulation should be related to
but but appears different from Kontsevich's homological mirror symmetry.
He used Fukaya's $A_{\infty}$ category of a symplectic manifold,
while we use an $A_{\infty}$ algebra structure on the de Rham cohomology
and its deformation constructed from the symplectic structure.
This has the advantage that on the mirror manifold,
which presumably is a complex manifold,
we has an $A_{\infty}$ algebra and its deformation constructed in the same
fashion.
This opens the door for establishing the identification 
by ideas from algebraic topology such as quasi-isomorphisms
(cf. Cao-Zhou \cite{Cao-Zho3, Cao-Zho4}).

The notion of an $A_{\infty}$ algebra was introduced by Stasheff \cite{Sta},
while that of an $L_{\infty}$ algebra first appeared in 
Schlessinger-Stasheff \cite{Sch-Sta}.
Such algebraic structures and their cousins show up in  various contexts in
string theory.
See Table 1 in Kimura-Stasheff-Voronov \cite{Kim-Sta-Vor} for
a guide to the literature. 
It is interesting to see whether the method of this paper yields 
the same infinite algebra structures as those 
in the literature by axiomatic approaches to conformal field theories.

We give a brief sketch of the history of homological perturbation theory here.
For more details,
see Gugenheim-Lambe-Stasheff \cite{Gug-Lam-Sta1}  
and the references therein.
One of Chen's major contributions was a method of computing the 
homology of the loop space on a manifold
by using his famous {\em formal power series connection}
and {\em iterated integrals} \cite{Che2}.
A purely algebraic version of Chen's result was
given by Gugenheim \cite{Gug},
and connection with $A_{\infty}$ algebra was maded 
by Gugenheim and Stasheff \cite{Gug-Sta}.
Then followed
much of the recent work on homological perturbation theory.
There were also independent development in the former Soviet Union.

This work is inspired by a recent manuscript by Huebschmann and 
Stasheff \cite{Hue-Sta} where they discussed the use of homological 
perturbation theory (HPT) in obtaining a formal power series connection on 
a DGLA.
There are some major difference between the content of this paper and
their work.
They focused on the existence of a formal power series connection
while we emphasize the deformations obtained from it.
Also they established the existence by basic perturbation lemma,
while we use the obstruction method following Hain \cite{Hai}.
We use basic perturbation lemma only to obtain the deformations
of the $A_{\infty}$ or $L_{\infty}$ structures.

The rest of the paper is arranged as follows.
It can be naturally divided into four parts.
In the first part (\S 1 - \S 5),
we recall the definition of formal power series connections
(\S 1),
the constructions given by Chen \cite{Che1} for DGA's (\S 2)
and by Hain \cite{Hai} for DGLA's (\S 3),
the relationships with $A_{\infty}$ algebras and $L_{\infty}$ algebras
respectively (\S 4)
and some basics of homological perturbation theory (\S 5).
In the second part (\S 6 - \S 9),
we apply the methods in the first part to various algebraic structures.
In the third part (\S 10 -\S 12),
we discuss some algebraic structures arising in various geometric contexts.
In the last part (\S 13),
we explain how extended deformation problems 
leads to the algebraic structures 
discussed in the third part,
then formulate the $A_{\infty}$ mirror symmetry conjecture 
by the applying the methods in the second part.

{\bf Acknowledgement}. 
{\em This work is partially supported by NSF group infra-structure
grant through the GAT group at Texas A\&M University.
I thank the Mathematics Department and members of the GAT group for 
hospitality.
Special thanks are due to Huai-Dong Cao for encouragement and collaboration
which leads to this work.
I also thank Jim Stasheff for sharing the manuscript \cite{Hue-Sta},
the reading of which greatly inspires me.}

\section{Formal power series connections}

Throughout this paper, $\bk$ is a graded commutative $\bQ$-algebra.

\subsection{Differential graded algebras}
Let $\cA = \sum_{n \in \bZ} \cA_n$ be a $\bZ$-graded module over $\bk$. 
If $\cdot: \cA \otimes \cA \to \cA$ is a $\bk$-linear map 
such that $\cA_m \cdot \cA_n \subset \cA_{m+n}$,
and $(a \cdot b) \cdot c = a \cdot (b \cdot c)$ for $a, b, c \in \cA$,
then $(\cA, \cdot)$ is called a {\em graded algebra}.
A module map $\phi: \cA \to \cA$ is called a {\em derivation of degree $k$}
if $\phi(\cA_n) \subset \cA_{n+k}$ for all $n \in \bZ$
and 
$$\phi(a \cdot b) = \phi(a) \cdot b + (-1)^{k |a|} a \cdot \phi(b),$$
for any homogeneous $a, b \in \cA$.
A {\em differential} is a derivation $\delta$ of degree $1$ such that 
$\delta^2 = 0$. 
A {\em differential graded algebra} (DGA) 
is a graded algebra $\cA$ together 
with a differential $\delta$,
and its {\em cohomology} is $H(\cA, \delta) = \Ker \delta/\Img \delta$.
The multiplication on $\cA$
induces a multiplication on $H(\cA, \delta)$.
It is easy to see that this multiplication is associative so
$H(\cA, \delta)$ has an induced structure of a graded algebra.

\subsection{Formal power series connections}
For a graded $\bk$-module $V$, 
let $\overline{T}(V) = \sum_{n \geq 0} T^n(V)$ 
be the completion of the graded tensor algebra on $V$.
For elements $v_1, \cdots, v_n \in V$,
we will write $v_1\cdots v_n$ for $v_1 \otimes \cdots \otimes v_n$.
A {\em formal power series connection} 
on a DGA $(\cA, \delta)$ is an element in $\cA \otimes \overline{T}(V)$
of the form
$$\omega = \sum \alpha_i v^i +
 \cdots + \sum \alpha_{i_1\cdots i_n} v_{i_1}\cdots v_{i_n} + \cdots,$$
where $v_{i_1}, \cdots, v_{i_n} \in V$,
$\alpha_{i_1\cdots i_n} \in \cA$
and $|\alpha_{i_1\cdots i_n}| = 1 - |v_{i_1}| - \cdots |v_{i_n}|$
(i.e. $\omega$ is of total degree $1$).
The multiplication on $\cA$ and 
the tensor product on $\overline{T}(V)$
induce a multiplication on $\cA \otimes \overline{T}(V)$
by the Koszul convention:
$$(\alpha v_1 \cdots v_m) \cdot (\beta w_1 \cdots w_n) 
= (-1)^{(|v_1| + \cdots + |v_m|)|\beta|} 
(\alpha \cdot \beta)(v_1 \cdots v_m w_1 \cdots w_n).$$ 
Also extend $\fd$ by
$$\fd(\alpha v_1 \cdots v_m) = (\fd\alpha) (v_1 \cdots v_m).$$
Then it is straightforward to check that $\fd$ is a derivation of 
$(\cA \otimes \overline{T}(V), \cdot)$.
We define the {\em curvature} of a connection $\omega$ to be
$$\Omega = \fd \omega + \omega \cdot \omega.$$
Formal power series connections and their curvatures were
introduced by Chen (see e.g. \cite{Che1}).
For our purpose,
we have used slightly a different definition of the curvatures.

\subsection{Some noncommutative formal deformations}
Given a connection $\omega$,
we consider the {\em covariant derivative} 
$\delta_{\omega}: \cA \otimes \overline{T}(V) 
\to \cA \otimes \overline{T}(V)$
defined by
$$\fd_{\omega} \alpha = \fd \alpha  + \omega \cdot \alpha$$
for $\alpha \in \cA \otimes \overline{T}(V)$.
Then we have
$$\fd_{\omega}^2 \alpha 
= (\fd + \omega \cdot)(\fd \alpha + \omega \cdot \alpha)
=(\fd \omega + \omega \cdot \omega) \cdot \alpha = \Omega \cdot \alpha. $$
When $\Omega = 0$,
$(\cA \otimes \overline{T}(V), \fd_{\omega})$ is a cochain complex.
Similarly,
if $\partial$ is a differential on $\overline{T}(V)$,
then we extend it to $\cA \otimes \overline{T}(V)$ by
$$\partial(\alpha v_1\cdots v_n ) 
= (-1)^{|\alpha|}\alpha \partial(v_1\cdots v_n),$$
where $\alpha \in \cA$, $v_1, \cdots, v_n \in V$.
Then $\partial$ is a differential on $\cA \otimes \overline{T}(V)$.
Since $\partial\fd + \fd \partial = 0$ on $\cA \otimes \overline{T}(V)$,
$(\partial + \fd)^2 = 0$
and so $(\cA \otimes \overline{T}(V), \partial + \fd)$ is a chain complex.
Given a connection $\omega$ 
and a differential $\partial$ on $\overline{T}(V)$,
$(\partial + \fd_{\omega})^2 
= \omega \partial + \fd \omega + \omega \cdot \omega$.
Hence if the latter vanishes,
$(\cA \otimes \overline{T}(V), \partial + \fd)_{\omega}$ 
is a chain complex.

\section{Chen's construction}
\label{sec:ChenConstruction}

\subsection{Cohomological splitting}
First fix a decomposition $\cA = \cH \oplus \fd M \oplus M$,
such that $\cH \subset \Ker \fd$, 
the natural map $\cH \to H(\cA, \fd)$ is
an isomorphism and $\fd|_M$ is injective.
Such a decomposition is called a {\em cohomological splitting}.
Since every element $a \in \cA$ can be uniquely written as
$a = a^H + \fd x_1 + x_2$ for some $a^H \in \cA$, $x_1, x_2 \in M$,
define
$Q: \cA \to \cA$ by $Q(a) = x_1$.
Then we have 
\begin{eqnarray*}
(1 - [\fd, Q])a = a - \fd Q a - Q\fd a = a - \fd x_1 - Q \fd x_2 =
a - \fd x_1 - x_2 = a^H,
\end{eqnarray*}
for all $a \in \cA$.
As a consequence,
if $\fd a = 0$,
then $a = a^H + \fd Q a$.

\subsection{Hodge theory and cohomological splitting}

Suppose that a DGA $(\cA, \cdot, \fd)$ is given
a (Euclidean or Hermitian) metric $\langle \cdot, \cdot \rangle$,
such that $\fd$ has a formal adjoint $\fd^*$,
i.e.,
$$\langle \fd a, b \rangle = \langle a, \fd^* b\rangle.$$
Since $\fd^2 = 0$,
it follows that $(\fd^*)^2 = 0$.
Set $\square_{\fd} = \fd\fd^* + \fd^* \fd$,
and $\cH = \Ker \square_{\fd} = \{a \in \cA: \square_{\fd} a = 0 \}$.
Then $\cH = \{a: \fd a = 0, \fd^* a = 0 \} = \Ker \fd \cap \Ker \fd^*$.
Assume that $A$ admits a ``Hodge decomposition'':
$$\cA = \cH \oplus \Img \fd \oplus \Img \fd^*.$$
It is standard to see that $\Ker \fd = \cH \oplus \Img \fd$,
and hence $H(\cA, \fd) \cong \cH$.
Note that $\square_{\fd}|_{\Img \fd \oplus \Img \fd^*}$ is invertible,
denote its inverse by $\square^{-1}$.
Define the Green's operator $G: \cA \to \cA$ as the composition
$$\cA \to \Img \fd \oplus \Img \fd^* 
\stackrel{\square^{-1}}{\longrightarrow}
\Img \fd \oplus \Img \fd^*\to \cA,$$
where the first arrow is the projection,
and the last arrow is the inclusion.
Since $\square_{\fd}$ commutes with $\fd$ and $\fd^*$,
so does $G$.
For any $\alpha \in \cA$,
denote its projection to $\cH$ by $\alpha^H$. 
Then we have 
\begin{eqnarray*}
\alpha - \alpha^H = \square_{\fd} G_{\fd} \alpha 
= (\fd\fd^* + \fd^*\fd) G_{\fd}\alpha
= \fd \fd^*G_{\fd}\alpha + \fd^* G_{\fd}\fd\alpha.
\end{eqnarray*}
This gives the Hodge decomposition explicitly.
Also if we set $M = \Img \fd^*$,
then $\Img \fd = \fd M$
and the Hodge decomposition is a cohomological splitting.
In this case, we have $Q = \fd^*Q_{\fd}$.

\subsection{Chen's construction}
Recall the following notations:
for a graded $\bk$-module $V$,
$sV$ is the graded $\bk$-module which has
$(sV)_n = V_{n+1}$ (raising the degree by $1$),
and $V^t$ the dual module of $V$ with its natural grading
(e.g. when $\bk$ is ungraded,
elements in $\Hom(W_n, \bk)$ has degree $-n$). 
Similarly define $s^{-1}$.
Then $(s^{-1}V)^t = sV^t$.

Chen \cite{Che1} inductively constructed a formal power series connection
$\omega \in \cA \otimes \overline{T}^n(s\cH^t)$.
Denote by $\omega_{[n]}$ the component of $\omega$ 
in $\cA \otimes T^n(s\cH^t)$,
$\partial_{[n]}$ the derivation determined 
by a map $s\cH^t \to T^n(s\cH^t)$.
Set $\omega_n = \omega_{[1]} + \cdots + \omega_{[n]}$
and $\partial_{[n]} = \partial_{[1]} + \cdots + \partial_{[n]}$.
Choose a homogeneous basis $\{\alpha_j \in \cH\}$,
denote by $\{X^j \in s\cH^t\}$ the dual basis of $\{s^{-1}\alpha_j\}$,
and $I$ the ideal in $\overline{T}(s\cH^t)$
or $\cA \otimes \overline{T}(s\cH^t)$ generated by $X^j$'s.
Take $\omega_{[1]} = \sum \alpha_j X^j$ and $\partial_{[1]} = 0$.
In general,
assume that $\partial_n$ and $\omega_n$ have been defined,
such that 
\begin{eqnarray*}
\partial_n \omega_n + \fd \omega_n + \omega_n \cdot \omega_n  = 0 
\pmod {I^{n+1}}, 
\end{eqnarray*}
set $\Gamma_{[n+1]} 
= \partial_n \omega_n + \fd \omega_n + \omega_n \cdot \omega_n
\pmod {I^{n+2}}$.
Then modifying a calculation in Chen \cite{Che1},
modulo $I^{n+2}$, we have 
\begin{eqnarray*}
\fd \Gamma_{[n+2]} 
& = & \fd (\partial_n \omega_n + \fd \omega_n + \omega_n \cdot \omega_n)
= \fd \partial_n \omega_n + \fd (\omega_n \cdot \omega_n) \\
& = & - \partial_n \fd \omega_n
+ \fd \omega_n \cdot \omega_n - \omega_n \cdot \fd \omega_n \\
& = & -\partial_n \fd \omega_n 
- (\partial_n \omega_n + \omega_n \cdot \omega_n) \cdot \omega_n 
+ \omega_n \cdot (\partial \omega_n + \omega_n \cdot \omega_n) \\
& = & -\partial_n \fd \omega_n - \partial_n \omega_n \cdot \omega_n 
+ \omega_n \cdot \partial_n \omega_n  
= - \partial_n (\fd \omega_n + \omega_n \cdot \omega_n) = 0.
\end{eqnarray*}
This shows $\Gamma_{[n+1]} \in \Ker \fd$,
so if one sets
\begin{align*}
& \partial_{[n+1]} \omega_{[1]} = - \Gamma_{[n+1]}^H,
& \omega_{[n+1]} & = - Q\Gamma_{[n+1]}, 
\end{align*}
then one has
\begin{eqnarray*}
\partial_{n+1} \omega_{n+1} 
+ \fd \omega_{n+1} + \omega_{n+1} \cdot \omega_{n+1}  = 0 \pmod {I^{n+2}}.
\end{eqnarray*}
Let $\partial = \partial_{[2]} + \cdots \partial_{[3]} + \cdots$
and $\omega = \Omega_1 + \omega_2 + \cdots$.

\begin{lemma} \label{lm:differential}
Let $n$ be a positive integer.
If $\partial$ is a derivation of $\overline{T}(s\cH^t)$ such that
$\partial I \subset I^2$,
and if $\omega \in \cA \otimes \overline{T}(s\cH^t)$ 
is a formal power seris connection such that 
\begin{itemize}
\item[(a)] $\omega \equiv \sum \alpha_j X^j \pmod I$,
\item[(b)] $\partial \omega + \fd \omega + \omega \cdot \omega 
	\equiv 0 \pmod {I^n}$,
\end{itemize}
then
\begin{itemize}
\item[(c)] $\fd (\partial \omega + \fd \omega + \omega \cdot \omega) 
	\equiv 0 \pmod {I^{n+1}}$,
\item[(d)] $\partial^2 \equiv 0 \pmod {I^{n+1}}$.
\end{itemize}
\end{lemma}

\begin{proof}
When $n = 1$,
the lemma is easily verified.
Suppose now $n> 1$ and the lemma holds for $n-1$.
Then we have
\begin{eqnarray*}
&& \fd (\partial \omega + \fd \omega + \omega \cdot \omega)
= \fd \partial \omega + \fd \omega \cdot \omega - \omega \cdot \fd \omega \\
& \equiv & - \partial \fd \omega 
- (\partial \omega + \omega \cdot \omega) \cdot \omega
+ \omega \cdot (\partial \omega + \omega \cdot \omega) \pmod {I^{n+1}} \\
& \equiv & - \partial (\fd \omega + \omega \cdot \omega) \pmod {I^{n+1}} 
	\equiv \partial^2 \omega \pmod {I^{n+1}}.
\end{eqnarray*}
Now by induction hypothesis,
$\partial^2 \equiv 0 \pmod {I^n}$,
this implies
$$\partial^2 \omega \equiv \sum \alpha_j (\partial^2 X^j) \pmod {I^{n+1}}.$$
Hence we have
$$ \fd (\partial \omega + \fd \omega + \omega \cdot \omega)
\equiv \sum \alpha_j (X^j \partial^2) \pmod {I^{n+1}}.$$
But the left hand side lies in $\Img \fd$ 
while the right hand side lies in $\cH$,
therefore,
they must both vanish modulo $I^{n+1}$.
\end{proof}

\begin{remark}
Chen's original proof of a version of this lemma (\cite{Che1}, Lemma 1.2.2)
used his famous transport formula and iterated integrals.
Gugenheim \cite{Gug} gave an algebraic proof.
Lemma 3.8 in Hain \cite{Hai} is the DGLA version of this lemma.
Here we use the DGA version of his proof.
\end{remark}

\section{Differential graded Lie algebra version}

One can also carry out the above discussion for DGLA's.
This has been done by Hain in his thesis under Chen's supervision 
(see e.g. \cite{Hai}).

\subsection{Differential graded Lie algebra}
Let $\cL = \sum \cL_n$ be a graded $\bk$-module,
a graded Lie algebra structure on $\cL$ is 
a degree zero map $[\cdot, \cdot]: \cL \otimes \cL \to \cL$,
such that 
\begin{eqnarray*}
&& [a, b] = -(-1)^{|a||b|} [b, a], \\
&& [a, [b, c]] = [[a, b], c] + (-1)^{|a||b|}[b, [a, c]],
\end{eqnarray*}
for homogeneous $a, b, c \in \cL$.
A module map $\phi: \cL \to \cL$ is called a {\em derivation of degree $k$}
if $\phi(\cL_n) \subset \cL_{n+k}$ for all $n \in \bZ$
and 
$$\phi([a, b]) = [\phi(a), b] + (-1)^{k |a|} [a, \phi(b)],$$
for any homogeneous $a, b \in \cL$.
In particular,
given any homogeneous $a \in \cL$,
$\ad_{a}: \cL \to \cL$ defined by $\ad_a (b) = [a, b]$
is a derivation of degree $|a|$.
A {\em differential} is a derivation $\fd$ of degree $1$
such that $\fd^2 = 0$,
a {differential graded Lie algebra} (DGLA) is 
a graded Lie algebra $(\cL, [\cdot, \cdot])$ together 
with a differential $\fd$.
The {\em cohomology} of a DGLA $(\cL, [\cdot, \cdot], \fd)$ is as usual
$H(\cL, \delta) = \Ker \fd /\Img \fd$.
The Lie bracket $[\cdot, \cdot]$ on $\cL$
induces a Lie bracket on $H(\cL, \delta)$
such that it becomes a graded Lie algebra.

\subsection{Formal power series connections on a DGLA}
Given a DGLA $(\cL, [\cdot, \cdot], \fd)$ and a graded $\bk$-module $V$, 
a {\em formal power series connection} 
is an element of total degree $1$ in $\cL \otimes \overline{S}(V)$
of the form
$$\omega = \sum \alpha_i v^i + \cdots 
+\sum \alpha_{i_1\cdots i_n} v_{i_1}\odot \cdots \odot v_{i_n} + \cdots,$$
where $v_{i_1}, \cdots, v_{i_n} \in V$
and $\alpha_{i_1\cdots i_n} \in \cL$.
The Lie bracket $[\cdot, \cdot]$ on $\cL$ and 
the symmetric product $\odot$ on $\overline{S}(V)$
induces a Lie bracket $[\cdot, \cdot]$ on $\cL \otimes \overline{S}(V)$ by:
$$[\alpha v_1 \odot \cdots \odot v_m, \beta w_1 \odot \cdots \odot w_n] 
= (-1)^{(|v_1| + \cdots + |v_m|)|\beta|} 
[\alpha, \beta]
(v_1 \odot \cdots \odot v_m \odot w_1 \odot \cdots \odot w_n).$$ 
Also extend $\fd$ by
$\fd(\alpha v_1 \cdots v_m) = (\fd\alpha) (v_1 \cdots v_m)$.
Then $(\cL \otimes \overline{S}(V), [\cdot, \cdot], \fd)$ is also a DGLA.
Given a connection $\omega$,
consider the {\em covariant derivative} 
$\delta_{\omega}: 
\cL \otimes \overline{S}(V) \to \cL \otimes \overline{S}(V)$
defined by
$$\fd_{\omega} \alpha = \fd \alpha  + [\omega, \alpha]$$
for $\alpha \in \cA \otimes \overline{S}(V)$.
Then we have
$$\fd_{\omega} [\alpha, \beta] 
= [\fd_{\omega} \alpha, \beta] 
+ (-1)^{|\alpha|} [\alpha, \fd_{\omega} \beta],$$
for homogeneous $\alpha, \beta \in \cL \otimes \overline{S}(V)$.
In other words,
$\fd_{\omega}$ is a derivation of degree $1$ on $\cL \otimes \overline{S}(V)$.
Furthermore,
$$\fd_{\omega}^2 \alpha 
= (\fd + \ad_{\omega})(\fd \alpha + [\omega,\alpha])
=[\fd \omega + \frac{1}{2}[\omega, \omega], \alpha] = [\Omega, \alpha],$$
where $\Omega = \fd \omega + \frac{1}{2}[\omega, \omega]$ 
is the {\em curvature}.
Hence if $\Omega = 0$,
$(\cL \otimes \overline{S}(V), \fd_{\omega})$ is a DGLA
which can be regarded as a formal deformation of $(\cL, \fd)$ as a DGLA.
Similarly,
if $\partial$ is a differential on $\overline{S}(V)$,
then we extend it to $\cL \otimes \overline{S}(V)$ by
$$\partial(\alpha v_1\odot\cdots\odot v_n ) 
= (-1)^{|\alpha|} \alpha \partial (v_1\odot\cdots\odot v_n),$$
where $\alpha \in \cL$, $v_1, \cdots, v_n \in V$.
Then 
$\partial$ is a differential on $\cL \otimes \overline{S}(V)$;
furthermore, $[\partial, \fd] =0$ on $\cL \otimes \overline{S}(V)$,
hence $(\partial + \fd)^2 = 0$.
If $\omega$ is a formal power series connection,
then one can consider the derivation $\partial + \fd_{\omega}$.
It is a differential if 
$\omega \partial + \fd \omega + \frac{1}{2}[\omega, \omega] = 0$.

Following Hian \cite{Hai},
given a cohomological splitting $\cL = \cH \oplus \fd M \oplus M$ of a 
DGLA $(\cL, [\cdot, \cdot], \fd)$,
the construction in \S \ref{sec:ChenConstruction} can be easily modified
to obtain a formal power series connection 
$\omega$ and a differential $\partial$ on $\overline{S}(s\cH^t)$ 
such that $\omega \partial + \fd \omega + \frac{1}{2}[\omega, \omega] = 0$.
Note in the proof of DGLA version of Lemma \ref{lm:differential},
one can use the fact that $[\omega, [\omega, \omega]] = 0$ 
for a power series connection.
There is no need to invoke the universal enveloping algebra 
as in Hain \cite{Hai}.

\section{Relationships with $A_{\infty}$ ($L_{\infty}$) 
algebras and twisting cochains}

\subsection{$A_{\infty}$ algebras}
An {\em $A_{\infty}$ algebra structure} on a graded $\bk$-module 
$V =\sum V_n$ is
a sequence of $\bk$-module maps
$m_n: \otimes^n V \to V$ of degree $2-n$,
such that for $v_1, \cdots, v_n \in V$,
\begin{eqnarray*} 
\sum_{r + s = n + 1} \sum_{k=1}^r(-1)^{\epsilon(k, s, v)} 
m_r(v_1, \cdots, m_s(v_k, \cdots, v_{k+s-1}), \cdots, v_n) = 0,
\end{eqnarray*}
where $\epsilon(k, s, v) = (s+1)k + s(n + |v_1| + \cdots + |v_{k-1}|)$.
The notion of an $A_{\infty}$ algebra was introduced in Stasheff \cite{Sta}.
It provides a generalization of the notion of a graded algebras and 
the notion of a DGA.
Indeed,
a graded algebra is equivalent to an $A_{\infty}$ algebra with 
$m_n = 0$ for $n \neq 2$,
a DGA is equivalent to an $A_{\infty}$ algebra with 
$m_n = 0$ for $n > 2$.

\subsection{Description of 
$A_{\infty}$ algebra structures by codifferentials}
There is a well-known equivalent description of $A_{\infty}$ structures
in terms of codifferentials on cotensor coalgebras
(the tilde construction in Stasheff \cite{Sta}).
Recall that a {\em graded coalgebra} is a graded $\bk$-module $C = \sum C_n$
together with a degree $0$ map
$\Delta: C \to C \otimes C$
such that the diagram
$$\CD
C @>{\Delta}>> C \otimes C \\
@V{\Delta}VV @VV{\Delta \otimes 1}V \\
C \otimes C @>>{1 \otimes \Delta}> C \otimes C \otimes C.
\endCD $$
commutes.
A {\em coderivation of degree $r$} on a coalgebra $(C, \Delta)$ is a map 
$L: C \to C$ of degree $r$ such that the diagram
$$\CD
C @>{L}>> C \\
@V{\Delta}VV @VV{\Delta}V \\
C \otimes C @>>{L \otimes 1 + 1 \otimes L}> C \otimes C
\endCD $$
commutes.
Here we use the Koszul convention,
e.g. $(1 \otimes L) (c_1 \otimes c_2) 
= (-1)^{|L||c_1|} c_1 \otimes L(c_2)$,
for $c_1, c_2 \in C$.
A {\em codifferential} is a coderivation $b$ of degree $1$ 
such that $b^2 = 0$.

For a graded $\bk$-module $W$,
the {\em cotensor coalgebra} on $W$ is the graded $\bk$-module
$T^c(W) = \sum W^{\otimes n}$ together with the comultiplication
$$\Delta (w_1 \otimes \cdots \otimes w_n) 
= \sum_{i=0}^n (w_1 \otimes \cdots \otimes w_i) \otimes 
(w_{i+1} \otimes \cdots \otimes w_n),$$
for $w_1, \cdots, w_n \in W$.
A coderivation $L: T^c(W) \to T^c(W)$
is uniquely determined by its composition with the projection
$\pi : T^c(W) \to W$.
Indeed,
denote by $L_{[n]}$ the restriction to $T^n(W)$ of $\pi\circ L$,
for simplicity, we assume $L_{[0]} = 0$, i.e., $L(1) = 0$,
then
\begin{eqnarray*}
&& L(w_1 \otimes  \cdots \otimes w_n) \\
& = & \sum_{i=1}^n \sum_{j = 1}^{n-i}
(-1)^{|L|(|w_1| + \cdots + |w_{j-1}|)}
w_1 \otimes \cdots \otimes 
L_{[i]}(w_j, \cdots, w_{i+j-1}) \otimes \cdots \otimes w_n,
\end{eqnarray*}
for $w_1, \cdots, w_n \in W$.
Furthermore,
if $L$ is a coderivation of odd degree with $L_{[0]} = 0$ on $T^c(W)$,
then $L^2$ is a coderivation of even degree,
and
\begin{eqnarray*}
&& (L^2)_{[n]}(w_1 \otimes \cdots \otimes w_n)  \\
& = & \sum_{r+s = n+1} \sum_{k = 1}^r (-1)^{|w_1| + \cdots + |w_{j-1}|}
L_{[r]}(w_1, \cdots, L_{[s]}(w_k, \cdots, w_{k+ s - 1}), \cdots,  w_n).
\end{eqnarray*}
Given a sequence of $\bk$-module maps $\{m_n: V^{\otimes n} \to V\}$ 
of degrees $2 -n$ on $V$,
define a coderivation $b$ of degree $1$ on $T^c(s^{-1}V)$ by
$$b_{[k]}(s^{-1}v_1, \cdots, s^{-1}v_k)
= (-1)^{\sum_{j=1}^k(k- j)|v_j|} s^{-1} m_k(v_1, \cdots, v_k).$$
(This is taken from Getzler-Jones \cite{Get-Jon},
Proposition 1.3. Unfortunately,
they included an extra $k(k-1)/2$ by mistake.)
Then $b^2 = 0$ if and only if
$\{ m_n \}$ is an $A_{\infty}$ algebra structure on $V$.

\subsection{Description of 
$A_{\infty}$ algebra structures by differentials}

Denote by $\langle \phi, w \rangle$ the pairing of an element $\phi \in W^t$ 
with an element $w \in W$.
Define the pairing between $(W^t)^{\otimes n}$ and $W^{\otimes n}$ by
\begin{eqnarray*}
\langle \phi^1 \otimes \cdots \otimes \phi^n, w_1 \otimes \cdots \otimes w_n
\rangle
= (-1)^{\sum_{j = 1}^{n-1} \sum_{k = j+1}^n |w_j||\phi^k|}
\langle \phi^1, w_1\rangle \cdots \langle \phi^n, w_n\rangle,
\end{eqnarray*}
where $\phi^1, \cdots, \phi^n \in W^*$,
$w_1, \cdots, w_n \in W$.
There is an one-to-one correspondence between 
$\End\overline{T}(W^t)$ and $\End\overline{T}(W)$ given by
$f \mapsto f^*$,
where $f \in \End \overline{T}(W)$, 
$f^* \in \End\overline{T}(W^t)$ is defined by
$$\langle f^*(\alpha) , X \rangle
= (-1)^{|f||\alpha|} \langle \alpha, f(X) \rangle,$$
where $\alpha \in \overline{T}(W^t)$ and $X \in \overline{T}(W)$.
It is straightforward to see that a derivation 
of degree $k$ on $\overline{T}(W)$
correspond to a derivation of degree $k$ on $\overline{T}(W^t)$.
Hence, there is a one-to-one correspondence between
the $A_{\infty}$ algebra structures on $V$ with differentials on
$\overline{T}((s^{-1}V)^t) = \overline{T}(sV^t)$.
In particular,
the differential $\partial$ on $\overline{T}(s\cH^t)$ constructed 
in \S \ref{sec:ChenConstruction} corresponds to an $A_{\infty}$
algebra structure on $\cH$ hence on $H(\cA, \fd)$.

\begin{example}
Given any closed connected oriented smooth manifold $X$,
the de Rham algebra is the DGA $(\Omega^*(X), \wedge, d)$. 
Its cohomology $H^*(X)$ is the de Rham cohomology of $X$.
When $X$ is endowed with a Riemannian metric and an orientation,
the classical Hodge theory provides the Hodge decomposition 
$$\Omega^*(X) = \cH \oplus \Img d \oplus \Img d^*.$$
Then by \S \ref{sec:ChenConstruction},
we get an $A_{\infty}$ algebra structure on $H^*(X)$.
\end{example}

\subsection{Twisting cochains}

If $C$ is a graded coalgebra with a differential $b$,
and $A$ is a DGA with differential $d$,
a {\em twisting cochain} is a $\bk$-module map $\tau: C \to A$ of degree $1$
such that $\tau b + d \tau = \tau \cup \tau$.
Here we have used the cup product 
$\cup: \Hom(C, A) \otimes \Hom(C, A) \to \Hom(C, A)$
defined by 
$$f \cup g = m(f \otimes g)\Delta,$$
for $f, g \in \Hom(C, A)$,
where $m$ denotes the multiplication on $A$.
It is clear that the formal power series connection constructed 
in \S \ref{sec:ChenConstruction}
corresponds to a twisting cochain 
$\tau: (T^c(s^{-1}H^*(\cA, \delta)), b) \to (\cA, \delta)$.

\subsection{$L_{\infty}$ algebras}
Recall that a $(k, n)$-shuffle is a permutation $\sigma$ 
which preserves the order of the sets 
$\{1, \cdots, k\}$ and $\{k+1, \cdots, n \}$
in the sense that whenever $\sigma(r) < \sigma(s) \leq k$ or 
$k+1 \leq \sigma(r) < \sigma(s)$,
then $r < s$.
A $(k, n)$-unshuffle is the inverse of an $(k, n)$-shuffle,
i.e. it is a permutation $\sigma$ of 
$\{1, 2, \cdots, n\}$ such that
$\sigma(1) < \cdots < \sigma(k)$ 
and $\sigma(k+1) < \cdots < \sigma(n)$.
Denote by $Sh^u(k, n)$ the set of all $(k, n)$-unshuffles.

An {\em $L_{\infty}$ algebra structure} on a graded $\bk$-module 
$V =\sum V_n$ is
a sequence of graded anti-symmetric $\bk$-module maps
$l_n: \otimes^n V \to V$ of degree $2-n$,
such that for $v_1, \cdots, v_n \in V$,
\begin{eqnarray*} 
\sum_{r + s = n + 1} \sum_{\sigma \in Sh^u(s, n)}
(-1)^{\sigma}(-1)^{\epsilon(\sigma, v)} 
l_r(l_s(v_{\sigma(1)}, \cdots, v_{\sigma(s)}), 
v_{\sigma(s+1)}, \cdots, v_{\sigma(n)}) = 0,
\end{eqnarray*}
where $(-1)^{\sigma}$ is the sign of the permutation $\sigma$,
and $\epsilon(\sigma, v)$ is determined by the Koszul convention.
The notion of an $L_{\infty}$ algebra first appeared in 
Schlesinger-Stasheff \cite{Sch-Sta}.
It provides a generalization of the notion of a graded Lie algebra 
and the notion of a DGLA.
Indeed,
a graded Lie algebra is equivalent to an $L_{\infty}$ algebra with 
$l_n = 0$ for $n \neq 2$,
a DGLA is equivalent to an $L_{\infty}$ algebra with 
$l_n = 0$ for $n > 2$.

\subsection{Description of 
$L_{\infty}$ algebra structures by codifferentials and differentials}
For a graded $\bk$-module $W$,
the {\em free cocommutative coalgebra} cogenerated by $W$ 
is the graded $\bk$-module
$\overline{S}^c(W) = \sum S^n(W)$ together with the comultiplication
\begin{eqnarray*}
&& \Delta^s (w_1 \odot \cdots \odot w_n) \\
& = & \sum_{i=0}^n \sum_{\sigma \in Sh^u(i, n)} 
(-1)^{\epsilon(\sigma, w)} 
(w_{\sigma(1)} \odot \cdots \odot w_{\sigma(i)}) \otimes 
(w_{\sigma(i+1)} \odot \cdots \odot w_{\sigma(n)}),
\end{eqnarray*}
for $w_1, \cdots, w_n \in W$.
A coderivation $L: \overline{S}^c(W) \to \overline{S}^c(W)$
is uniquely determined by its composition with the projection
$\pi : \overline{S}^c(W) \to W$.
Indeed,
denote by $L_{[n]}$ the restriction to $S^n(W)$ of $\pi\circ L$,
for simplicity, we assume $L_{[0]} = 0$, i.e., $L(1) = 0$,
then
\begin{eqnarray*}
&& L(w_1 \odot  \cdots \odot w_n) \\
& = & \sum_{s=1}^n \sum_{\sigma \in Sh^u(s, n)}
(-1)^{\epsilon(\sigma, w)}
L_{[i]}(w_{\sigma(1)}, \cdots, w_{\sigma(i)}) \odot 
w_{\sigma(i+1)} \odot \cdots \odot w_{\sigma(n)},
\end{eqnarray*}
for $w_1, \cdots, w_n \in W$.
Furthermore,
if $L$ is a coderivation of odd degree with $L_{[0]} = 0$ on $T^c(W)$,
then $L^2$ is a coderivation of even degree,
and
\begin{eqnarray*}
&& (L^2)_{[n]}(w_1 \odot \cdots \odot w_n)  
= \sum_{r+s = n+1} \sum_{\sigma \in Sh^u(s, n)}
(-1)^{\epsilon(\sigma, w)} \\
&& \qquad L_{[r]}(L_{[s]}(w_{\sigma(1)}, \cdots, w_{\sigma(s)}), 
w_{\sigma(s+1)}, \cdots, \cdots,  w_{\sigma(n)}).
\end{eqnarray*}
Given a sequence of graded anti-symmetric $\bk$-module maps 
$\{l_n: V^{\otimes n} \to V\}$ of degrees $2 - n$ on $V$,
each $l_n$ can be regarded as a map of degree $1$ 
from $(s^{-1}V)^{\otimes n}$ to $s^{-1}V$,
where as usual $s^{-1}V$ is the graded $\bk$-module which has
$(s^{-1}V)_n = V_{n+1}$.
Define a coderivation $b$ of degree $1$ on $\overline{S}^c(s^{-1}V)$ by
$$b_{[k]}(s^{-1}v_1, \cdots, s^{-1}v_k)
= (-1)^{\sum_{j=1}^k(k- j)|v_j|} s^{-1} l_k(v_1, \cdots, v_k).$$
Then each map $b_{[k]}$ is symmetric on $(s^{-1}V)^{\otimes n}$ 
(cf. Penkava \cite{Pen}).
Furthermore,
$b^2 = 0$ if and only if
$\{ l_n \}$ is an $L_{\infty}$ algebra structure on $V$.
Dualizing,
one sees that an $L_{\infty}$ algebra structure on $V$
corresponds to a differential 
on $\overline{S}((s^{-1}V)^t)=\overline{S}(sV^t)$.

\section{Homological perturbation theory}

Chen's work inspired the obstruction method 
in homological perturbation theory.
Another important inspiration 
is the basic perturbation lemma
which first occurred in Brown \cite{Bro}
and independently in Gugenheim \cite{Gug}.
Both authors were inspired by Shih \cite{Shi}.

\subsection{SDR data}
Suppose that we are given two  chain complexes $(M, d_M)$ and $(A, d_A)$,
chain maps
$\nabla: M \to A$, $f: A \to M$ and a chain homotopy $\phi: A \to A$
such that 
\begin{eqnarray}
& f \nabla = id_M, \label{SDR1} \\ 
& \nabla f = id_A + d_A \phi + \phi d_A.\label{SDR2}
\end{eqnarray}
All this information is called an {\em SDR data},
and denoted by
$$\left(M \underset{f}{\overset{\nabla}\rightleftarrows} A, \phi\right).$$
Given an SDR data, 
one can adjust $\phi$ such that the following additional conditions 
are also satisfied (Lambe and Stasheff \cite{Lam-Sta}):
\begin{align*}
& \phi \nabla = 0, 
& f \phi & = 0, 
& \phi^2 & = 0.
\end{align*}
For example,
suppose that we have a cochain complex $(\cA, \fd)$ and 
a cohomological splitting $\cA = \cH \oplus \fd M \oplus M$,
then we take $A = \cA$, $M = \cH$,
$\nabla$ the inclusion,
$f$ the projection and $\phi = -Q$.
This gives an SDR which satisfies the additional conditions.
 We will need the following:

\begin{theorem} (Gugenheim and Stasheff \cite{Gug-Sta})
Given an SDR data
$$\left(M \underset{f}{\overset{\nabla}\rightleftarrows} A, 
\phi\right)$$
with $A$ a DGA,
there exist two maps of degree $-1$,
$\partial: T^c(sM) \to T^c(sM)$ 
and $\tau: T^c(sM) \to A$,
such that $\partial$ is a coderivation on $T^c(M)$ such that $\partial^2 = 0$,
and $\tau$ is a twisting cochain,
i.e. $\tau \partial + \fd \tau = \tau \cup \tau$.
\end{theorem}

This can be proved by the obstruction method,
but a delicate filtration is required to carry out the induction.

\subsection{Basic perturbation lemma}
The setup for basic perturbation lemma is as follows.
Given SDR-data
$$\left((M, d_M) \underset{f}{\overset{\nabla}\rightleftarrows} (N, d_N), 
\phi\right)$$ 
and a new differential $\cD_N$ on $N$,
one seeks for a new SDR-data
$$\left((M, \cD_M)
\underset{f_{\infty}}{\overset{\nabla_{\infty}}\rightleftarrows} (N, \cD_N),
\phi_{\infty} \right)$$
Again, this is achieved inductively. 
Let $t = \cD_N - \fd_N$, and call it the {\em initiator}.
For $n \geq 0$,
construct sequences $\{\nabla_n\}$, $\{f_n\}$, $\{\cD_{M_n}\}$ and $\{\phi_n\}$
inductively by
\begin{align*}
\nabla_{n+1} & = r_n \nabla, &&  f_{n+1} = f s_n, \\
\cD_{M_n} & = \fd_M + f \Sigma_n \nabla, &&  \phi_n = \phi s_n,
\end{align*}
where
\begin{align*}
& t_{n+1} = (t\phi)^n t, &&  \Sigma_n = t_1 + \cdots + t_n, \\
& s_n = 1 + \Sigma_n \phi, &&  r_n = 1 + \phi \Sigma_n.
\end{align*}
Assume that there are complete filtrations on $M$ and $N$ 
such that $\nabla$, $f$, $\partial$ and $\phi$ preserve the filtrations
and $t$ lowers the filtration,
then the relevant sequences converge and one obtains a new SDR.
A surprising fact is that the above construction actually preserves the 
algebra or coalgebra structure.
See Gugenheim, Lambe and Stasheff \cite{Gug-Lam-Sta1, Gug-Lam-Sta2}
for details.

\section{Applications to DGA's and DGLA's}
\label{sec:DGADGLA}

Given a DGA $(\cA, \cdot, \fd)$ with a cohomological splitting
$\cA = \cH \oplus \fd M \oplus M$,
one gets an SDR-data
\begin{eqnarray} \label{eqn:SDR1}
\left((\cH, 0) \underset{f}{\overset{\nabla}\rightleftarrows} 
(\cA, \fd), \phi\right),
\end{eqnarray}
with $\nabla$ the inclusion, 
$f$ the projection and $\phi = - Q$.
Tensoring with $\overline{T}(s\cH^t)$,
one gets a new SDR-data
\begin{eqnarray} \label{eqn:SDR2}
\left((\cH \otimes \overline{T}(s\cH^t), 0) 
\underset{f}{\overset{\nabla}\rightleftarrows} 
(\cA \otimes \overline{T}(s\cH^t), \fd), \phi\right)
\end{eqnarray}
From the algebra structure on $\cA$,
one obtains a derivation $\partial^a$ on $\overline{T}(s\cH^t)$ 
and a connection $\omega^a \in \cA \otimes \overline{T}(s\cH^t)$
such that $\partial^a \omega^a + \fd \omega^a + \omega^a \cdot \omega^a = 0$.
Consider the perturbation of $\fd$ by $\fd + \partial$,
we regard it as a noncommutative formal deformation of $\fd$,
then by basic perturbation lemma,
one gets a new SDR-data
\begin{eqnarray*} 
\left((\cH \otimes \overline{T}(s\cH^t), \cD) 
\underset{f_{\infty}}{\overset{\nabla_{\infty}}\rightleftarrows} 
(\cA \otimes \overline{T}(s\cH^t), \fd + \partial), \phi_{\infty} \right)
\end{eqnarray*}
Since on the right we have a DGA,
one can repeat the process.
This gives us a derivation $\partial^{aa}$ on 
$\overline{T}(s(\cH \otimes \overline{T}(s\cH^t))^t)$,
which we regard as a noncommutative formal deformation of $\partial^a$.
In other words,
$\partial^a$ gives an $A_{\infty}$ algebra structure on $\cH$,
$\partial^{aa}$ gives a noncommutative formal deformation of this structure.
One can repeat this process for as many times as one wishes.
One can also can consider the perturbations 
$\fd_{\omega^a}$ and $\partial^a + \fd_{\omega^a}$,
however,
these are not derivations.

One can do the same things for a DGLA to obtain (graded) formal deformations.
Also one can use $\fd_{\omega^L}$ and $\partial^L + \fd_{\omega^L}$
since they are also derivations.

\section{Differential Poisson algebra}
\label{sec:DiffPA}

Recall that a {\em graded Poisson algebra structure} on a graded $\bk$-module
consists of a multiplication $\cdot: \cA \otimes \cA$
and a Lie bracket $[\cdot, \cdot]: \cA \otimes \cA \to \cA$,
such that $(\cA, \cdot)$ is a graded associative algebra 
(we do not assume that the multiplication is graded commutative),
$(\cA, [\cdot, \cdot])$ is a graded Lie algebra,
and for homogeneous $a, b, c \in \cA$,
we have
$$[a, b \cdot c] = [a, b] \cdot c + (-1)^{|a||b|} b \cdot [a, c],$$
i.e.,
for any homogeneous $a \in \cA$,
$\ad_a = [a, \cdot]: \cA \to \cA$ is a derivation of degree $|a|$ of
$(\cA, \cdot)$.
A derivation of degree $k$ on a Poisson algebra $(\cA, \cdot, [\cdot, \cdot])$
is a map $\fd: \cA \to \cA$ of degree $k$ such that it is a derivation
for both $(\cA, \cdot)$ and $(\cA, [\cdot, \cdot])$.
As usual, a differential is a derivation $\fd$ of degree $1$ such that 
$\fd^2 = 0$.
A Poisson algebra with a differential is called a 
{\em differential Poisson algebra}.

There are variants to the discussions in \S \ref{sec:DGADGLA}
for a differential Poisson algebra.
Suppose that we are given a differential 
Poisson algebra $(\cA, \cdot, [\cdot, \cdot], \fd)$
and a cohomological splitting $\cA = \cH \oplus \fd M \oplus M$,
or equivalently,
we have an SDR-data
\begin{eqnarray} \label{eqn:SDR4}
\left((\cH, 0) \underset{f}{\overset{\nabla}\rightleftarrows} 
(\cA, \fd), \phi\right),
\end{eqnarray}
with $\nabla$ the inclusion, 
$f$ the projection and $\phi = - Q$.
Then using the associative algebra structure on $\cA$,
we obtain a differential $\partial^a$ on $\overline{T}(s\cH^t)$ and
a connection $\omega^a \in \cA \otimes \overline{T}(s\cH^t)$
such that $\partial^a \omega^a + \fd \omega^a + \omega^a \cdot \omega^a = 0$;
using the Lie algebra structure on $\cA$,
we obtain a differential $\partial^L$ on $\overline{S}(s\cH^t)$ and 
a connection $\omega^L \in \cA \otimes \overline{S}(s\cH^t)$
such that 
$\omega^L \partial^L + \fd \omega^L + \frac{1}{2}[\omega^L, \omega^L] = 0$.
Tensoring with $\overline{S}(s\cH^t)$,
we get from (\ref{eqn:SDR1}) an SDR-data
\begin{eqnarray} \label{eqn:SDR5}
\left((\cH \otimes \overline{S}(s\cH^t), 0) 
\underset{f}{\overset{\nabla}\rightleftarrows} 
(\cA \otimes \overline{S}(s\cH^t), \fd), \phi\right)
\end{eqnarray}
Now $\fd + \partial^L + \ad_{\omega^L}$ is a 
differential for the multiplicative structure on 
$\cA \otimes \overline{S}(s\cH^t)$,
it is a formal deformation of the differential $\fd$.
With $t = \partial^L + \ad_{\omega^L}$ as an initiator,
the basic perturbation lemma yields a new SDR-data
\begin{eqnarray*}
\left((\cH \otimes \overline{S}(s\cH^t), \cD) 
\underset{f_{\infty}}{\overset{\nabla_{\infty}}\rightleftarrows} 
(\cA \otimes \overline{S}(s\cH^t), \fd), \phi_{\infty} \right)
\end{eqnarray*}
which gives a differential $\partial^{aL}$ on 
$\overline{T}((s(\cH \otimes \overline{S}(s\cH^t)))^*)$.
In other words,
$\partial^a$ gives an $A_{\infty}$ algebra structure on $\cH$,
$\partial^{aL}$ gives a formal deformation of 
this $A_{\infty}$ algebra structure.
Similarly,
since $\fd + \partial^L$ is also a derivation,
we can also use $\partial^L$ as the initiator.

\section{Differential Gerstenhaber algebra}
\label{sec:DiffGA}

Recall that a {\em Gerstenhaber algebra} (or simply G-algebra) structure
on a graded $\bk$-module $\cA$
consists of an associative multiplication $\cdot$
(again we do not assume that the multiplication is graded commutative),
and a Lie algebra structure $[\cdot \bullet \cdot]$ on $s\cA$,
such that for homogeneous $a, b, c \in \cA$,
we have
$$[a \bullet (b \cdot c)] 
= [a \bullet b] \cdot c + (-1)^{(|a|-1)|b|} b \cdot [a \bullet c],$$
i.e.,
for any homogeneous $a \in \cA$,
$\ad_a = [a\bullet \cdot]: \cA \to \cA$ is a derivation of degree $|a| -1 $ of
$(\cA, \cdot)$.
A derivation of degree $k$ on a Gerstenhaber algebra 
$(\cA, \cdot, [\cdot \bullet \cdot])$
is a map $\fd: \cA \to \cA$ of degree $k$ such that it is a derivation
for both $(\cA, \cdot)$ and $(s\cA, [\cdot\bullet \cdot])$.
As usual, a differential is a derivation $\fd$ of degree $1$ such that 
$\fd^2 = 0$.
A Gerstenhaber algebra with a differential is called a 
{\em differential Gerstenhaber algebra}.

Suppose that we are given a differential 
Gerstenhaber algebra $(\cA, \cdot, [\cdot \bullet \cdot], \fd)$
and a cohomological splitting
$\cA = \cH \oplus \fd M \oplus M$,
we can proceed as in \S \ref{sec:DiffPA}.
We have an SDR-data
\begin{eqnarray} \label{eqn:SDR6}
\left((\cH, 0) \underset{f}{\overset{\nabla}\rightleftarrows} 
(\cA, \fd), \phi\right),
\end{eqnarray}
with $\nabla$ the inclusion, 
$f$ the projection and $\phi = - Q$.
Using the associative algebra structure on $\cA$,
we obtain a differential $\partial^a$ on $\overline{T}(s\cH^t)$;
using the Lie algebra structure on $s\cA$,
we obtain a differential $\partial^L$ on $\overline{S}(s\cH^t)$ and 
a connection $\omega^L \in \cA \otimes \overline{S}(s\cH^t)$
such that 
$\omega^L \partial^L + \fd \omega^L + \frac{1}{2}[\omega^L, \omega^L] = 0$.
Tensoring with $\overline{S}(s\cH^t)$,
we get from (\ref{eqn:SDR1}) an SDR-data
\begin{eqnarray} \label{eqn:SDR7}
\left((\cH \otimes \overline{S}(s\cH^t), 0) 
\underset{f}{\overset{\nabla}\rightleftarrows} 
(\cA \otimes \overline{S}(s\cH^t), \fd), \phi\right)
\end{eqnarray}
Now $\fd + \partial^L + \ad_{\omega^L}$ is a 
differential for the multiplicative structure on 
$\cA \otimes \overline{S}(s\cH^t)$,
it is a formal deformation of the differential $\fd$.
With $t = \partial^L + \ad_{\omega^L}$ as an initiator,
the basic perturbation lemma yields a new SDR-data
\begin{eqnarray*}
\left((\cH \otimes \overline{S}(s\cH^t), \cD) 
\underset{f_{\infty}}{\overset{\nabla_{\infty}}\rightleftarrows} 
(\cA \otimes \overline{S}(s\cH^t), \fd), \phi_{\infty} \right)
\end{eqnarray*}
which gives a differential $\partial^{aL}$ on 
$\overline{T}((s(\cH \otimes \overline{S}(s\cH^t)))^*)$.
In other words,
$\partial^a$ gives an $A_{\infty}$ algebra structure on $\cH$,
$\partial^{aL}$ gives a formal deformation of 
this $A_{\infty}$ algebra structure.

\section{(Differential) Gerstenhaber-Batalin-Vilkovisky algebras}

One interesting way to obtain G-algebras is to start 
with a graded commutative algebra with unit $(\cA, \wedge)$
and an odd operator $\Delta$ on it such that $\Delta^2 = 0$.
Set
\begin{eqnarray} \label{eqn:Tian}
[a \bullet b]_{\Delta} = (-1)^{|a|} (\Delta(a \wedge b) - \Delta a \wedge b
- (-1)^{|a|} a \wedge \Delta b),
\end{eqnarray}
for homogeneous $a, b \in \cA$.
It can be shown that if 
\begin{eqnarray*}
[a \bullet (b \wedge c)]_{\Delta} 
= [a \bullet b]_{\Delta} \wedge c 
+ (-1)^{(|a|-1)|b|}b \wedge [a \bullet c]_{\Delta},
\end{eqnarray*}
for homogeneous $a, b, c \in \cA$,
then one has
\begin{eqnarray*}
[a \bullet [b \bullet c]_{\Delta}]_{\Delta}
=[[a \bullet b]_{\Delta} \bullet c]_{\Delta}
+ (-1)^{(|a|-1)(|b|-1)} [b\bullet [a \bullet c]_{\Delta}]_{\Delta},
\end{eqnarray*}
i.e., $(\cA, \wedge, [\cdot \bullet \cdot]_{\Delta})$ is a G-algebra.
The tuple $(\cA, \wedge, \Delta, [\cdot\bullet\cdot]_{\Delta})$
is called a {\em Gerstenhaber-Batalin-Vilkovisky (GBV) algebra}.
Given such an algebra,
one easily sees that 
$$\Delta[a \bullet b]_{\Delta} = [\Delta a \bullet b]_{\Delta}
+ (-1)^{|a|-1} [a \bullet \Delta b]_{\Delta},$$
i.e. $(\cA[1], \Delta, [\cdot \bullet \cdot]_{\Delta})$ is a DGLA.
However, 
by (\ref{eqn:Tian}),
$[\cdot\bullet\cdot]_{\Delta}$ induces 
the trivial bracket on $H(\cA, \Delta)$.

A {\em differential Gerstenhaber-Batalin-Vilkovisky (DGBV) algebra}
is a GBV algebra $(\cA, \wedge, \Delta, [\cdot\bullet\cdot]_{\Delta})$ 
together with a derivation $\fd$ of odd degree with respect to 
$\wedge$,
such that $\fd^2 = [\fd, \Delta] = 0$.
It then follows that $\fd$ is also a derivation with respect to 
$[\cdot\bullet\cdot]_{\Delta}$ such that 
$(\cA, \wedge, \fd, [\cdot\bullet\cdot]_{\Delta})$ is a 
differential G-algebra.
Then one can apply the constructions in \S \ref{sec:DiffGA}.

\section{GBV algebra structure on the space of polyvector fields}
\label{sec:Poly}

Denote by $\Omega^{-*}(X) = \Gamma(X, \Lambda^*TX)$ 
the space of polyvector fields on $X$.
As $\Omega^*(X)$,
$\Omega^{-*}(X)$ has a wedge product $\wedge$ which makes it a
graded commutative associative algebra with unit.
The {\em Schouten-Nijenhuis bracket}
$[\cdot, \cdot]_S: \Omega^{-*}(X) \times \Omega^{-*}(X) 
	\to \Omega^{-*}(X)$
is characterized by the following properties:
\begin{itemize}
\item[(i)] for all $V \in \Omega^{-*}(X)$,
the endomorphism $\ad_V: u \to [V, u]_S$ is 
the Lie derivative of the vector field $V$,
\item[(ii)] if $u \in \Omega^{-p}(X)$ and $v \in \Omega^{-q}(X)$,
then we have
$$[u, v]_S = - (-1)^{(p-1)(q-1)}[v, u]_S.$$
\item[(iii)] if $u \in \Omega^{-*}(X)$,
$\ad_u$ is a derivation of degree $p-1$ of 
the algebra $(\Omega^{-*}(X), \wedge)$.
\end{itemize}
There are two ways to obtain DGBV algebra structures on $\Omega^{-*}(M)$.
Given any torsion free connection $\nabla$ on $TX$,
Koszul \cite{Kos} defined a generalized divergence operator  
$\Delta = - \sum_{i=1}^n \iota_{e^i} \nabla_{e_i}$,
where $\{e_1, \cdots, e_n\}$ is a local frame of $TX$,
and $\{e^1, \cdots, e^n\}$ is the dual frame.
He proved
\begin{eqnarray} 
&& \Delta^2 = 0, \label{eqn:differential} \\
&& [u, v]_S 
= (-1)^{|u|} (\Delta (u \wedge v) 
- \Delta u \wedge v - (-1)^{|u|} u \wedge \Delta v), \label{eqn:BV}
\end{eqnarray}
for homogeneous $u, v \in \Omega^{-*}(X)$.
I.e., 
$(\Omega^{-*}(X), \wedge, \Delta, [\cdot, \cdot]_S)$ is a GBV algebra.
Witten \cite{Wit} observed the following:
given an volume element $\vol \in \Omega^n(M)$,
let $\Delta_{\vol}$ be the conjugation of $d$ on $\Omega^{n-*}(M)$ by the 
isomorphism $\Omega^{-*}(M) \to \Omega^{n-*}(M)$ induced from $\vol$,
then $\Delta_{\vol}$ satisifes (\ref{eqn:differential}) 
and (\ref{eqn:BV}),
hence inducing a GBV algebra structure.
If one takes a torsion free connection which preserves the volume form,
then it is easy to see that $\Delta_{\vol} = \Delta$ for this
connection.

Let $g$ be a Riemannian metric on $X$,
then the Levi-Civita connection is torsion free.
The metric give an isomorphism $\flat: \Omega^{-*}(X) \cong \Omega^*(X)$
(and its inverse $\sharp$).
If $\{e_i\}$ is a local orthonormal frame of $TX$
and $\{e^i\}$ is the dual frame,
then $\flat$ is obtained by extending $e_i \mapsto e^i$.
Furthermore, if $\nabla$ is the Levi-Civita connection,
from the well-know formulas  
\begin{align*}
d & = \sum_{i=1}^n e^i \wedge \nabla_{e_i}, 
& d^* & = - \sum_{i=1}^n \iota_{e_j} \nabla_{e_j},
\end{align*}
we get $\Delta u = (d^* u^{\flat})^{\sharp}$.
When $M$ is endowed with an orientation,
the Riemannian metric defined a volume form $d \vol_g$.
Using this volume form,
one can define the isomorphism 
$\flat_{\vol}: \Omega^{-*}(X) \to \Omega^{n-*}(X)$ 
given by the contraction with $d \vol_g$.
Denote its inverse by $\sharp_{\vol}$,
then it is straightforward to see that 
$\Delta u = (d u^{\flat_{\vol}})^{\sharp_{\vol}}$.

\section{DGBV algebras from Poisson and symplectic manifolds}
\label{sec:Poisson}

\subsection{DGBV algebras from Poisson manifolds}

A {\em Poisson structure} is an elment $w \in \Omega^{-2}(X)$
such that $[w, w]_S = 0$.
Given a Poisson structure $w$,
$\ad_w$ is a differential on $\Omega^{-*}(X)$
and $[\Delta, \ad_w] = \ad_{\Delta w}$.
Therefore,
if $\Delta w = 0$,
then $(\Omega^{-*}(X), \wedge, \ad_w, \Delta, [\cdot,\cdot]_S)$
is a DGBV algebra.
This condition can be satisfied when $w$ is {\em regular},
i.e.,
when regarded as a family of anti-symmetric matrices on $\bR^n$,
$w$ has constant rank.
It is well-known that a regular Poisson manifold admits a torsion free
connection $\nabla$ such that $\nabla w = 0$.
Such a connection is called a {\em Poisson connection}.

Given a Riemannian metric on $X$,
one can use the isomorphism $\flat: \Omega^{-*}(X) \to \Omega^*(X)$
to transform the DGBV algebra structure on $\Omega^{-*}(X)$
to a DGBV algebra structure on $\Omega^*(X)$.
If the torsion free connection $\nabla$ is the Levi-Civita connection,
then $\Delta$ corresponds to $d^*$ on $\Omega^*(X)$.
This shows that
$(\Omega^*(M), \wedge, d^*, [\cdot \bullet \cdot]_{d^*})$ 
is a GBV algebra.
One can prove this by two other method,
one is by the local expression as in Cao-Zhou \cite{Cao-Zho1, Cao-Zho2},
the other is to notice that if $X$ is oriented,
the Riemannian metric determines a volume form $\vol_g$,
then $\Delta_{\vol_g}$ is conjugate with $d^*$ by $\flat$.

There is another important construction of DGBV algebra structure 
on $\Omega^*(X)$ of a Poisson manifold $(X, w)$.
Let $\Delta = [\iota_w, d]$ (Koszul \cite{Kos}),  
then 
$(\Omega^*(X), \wedge, d, \Delta, [\cdot\bullet\cdot]_{\Delta})$
is a DGBV algebra.
This DGBV algebra has a very special property.
If $\alpha, \beta \in \Omega^*(X)$ are closed,
then 
\begin{eqnarray*}
[\alpha\bullet\beta]_{\Delta} 
= (-1)^{|\alpha|-1}d (\iota_w(\alpha \wedge \beta)
- \iota_w \alpha \wedge \beta - \alpha \wedge \iota_w \beta),
\end{eqnarray*}
hence $[\cdot\bullet\cdot]_{\Delta}$ induces the zero bracket on $H^*(M)$.

\subsection{DGBV algebras from symplectic manifolds}
A symplectic structure on $X$ is a closed $2$-form $\omega$ nondegenerate 
in the sense that $\omega$ induces an isomorphism 
$\flat_{\omega}: T_xX  \to T_x^*X$ 
defined by $V \mapsto \omega(V, \cdot)$, for any $x \in X$, $V \in T_xX$,
with inverse $\flat_{\omega}$.
Extend $\flat_{\omega}$ and $\sharp_{\omega}$ on isomorphisms
between $\Omega^*(X)$ and $\Omega^{-*}(X)$.
Then $\omega^{\sharp_{\omega}} \in \Omega^{-2}(X)$ is a Poisson structure.
One can verify this by Darboux theorem.
As a Poisson manifold,
a symplectic manifold has the DGBV algebra structures 
on $\Omega^{-*}(X)$ and $\Omega^*(X)$ discussed above.
A special feature is that $d \alpha = [\omega \bullet \alpha]_{\Delta}$.

\section{DGBV algebras from complex geometry}

\subsection{Differential G-algebras from complex manifolds}

Given a complex $n$-manifold $M$,
denote by $T_cX$ the complex tangent bundle.
We have the following decomposition:
$$\Omega^{-*}(X) \otimes \bC
=\Omega^{-*, -*}(X) = \sum_{p, q\geq 0} \Omega^{-p, -q}(X) 
= \sum_{p, q\geq 0}  
\Gamma(M, \Lambda^pT_cX \otimes \Lambda^q\overline{T}_cX).$$
It is a complex G-algebra.
From $[\Omega^{-1, 0}(X), \Omega^{-1, 0}(X)] \subset \Omega^{-1, 0}(X)$,
one sees that $\Omega^{-*, 0}(X)$ is G-subalgebra of
$\Omega^{-*, -*}(X)$.
Using local complex coordinates,
it is easy to see that the Schouten-Nijenhuis bracket on
$\Omega^{-*, 0}(X)$ can be extended to a G-algebra structure on 
$\Omega^{-*, *}(X) 
= \Gamma(X, \Lambda^*T_cX \otimes \Lambda^*\overline{T}_c^*X)$.
Since  $\Omega^{-*, *}(X)$ is the space of $(0, *)$-forms on
the holomorphic vector bundle $\Lambda^* T_cX$,
there is a naturally defined differential $\bar{\partial}$ on it.
This makes $\Omega^{-*, *}(X)$ a differential G-algebra.

\subsection{DGBV algebras from Calabi-Yau manifolds}
If $\Omega$ is a nonvanishing section of $\Omega^{-n, 0}(X)$,
it defines an isomorphism $\Omega^{-*, *}(X) \to \Omega^{n-*, *}(X)$.
Let $\Delta: \Omega^{-*, *}(X) \to \Omega^{-(*-1), *}(X)$
be the conjugation of $\partial$ by this isomorphism,
then $\Delta$ gives a GBV algebra structure on $\Omega^{-*, *}(X)$.
If $\bar{\partial}\Omega = 0$,
then we have $\bar{\partial} \Delta + \Delta \bar{\partial} = 0$.
Hence 
$(\Omega^{-*, *}(M), \wedge, \bar{\partial}, \Delta, [\cdot, \cdot])$
is a DGBV algebra if $\Omega$ is a holomorphic volume form.
Recall that a Calabi-Yau manifold admits 
a parallel holomorphic volume form.
Such a construction of the operator $\Delta$ appeared 
in Tian \cite{Tia} and Todorov \cite{Tod} 
in the study of deformations of Calabi-Yau manifolds.
The relationship of $\Delta$ with the Schouten-Nijenhuis bracket on
$\Omega^{-*, *}(M)$ was revealed by a formula Tian proved.
These results actually appeared 
before Witten's observation in the real case \cite{Wit} 
which we mentioned in \S \ref{sec:Poly}.

\subsection{Differential Poisson 
algebras from holomorphic vector bundles}

Let $\pi: E \to X$ be a holomorphic vector bundle 
on a complex manifold $X$.
Since each fiber of $\End E$ is a (noncommutative) Poisson algebra,
one gets a graded Poisson algebra structure 
on $\Omega^{0, *}(X, \End E)$,
while the differential $\bar{\partial}_{\End E}$ makes it a 
differential Poisson algebra. 
Similarly for $\Omega^{*, *}(X, \End E)$
and $\Omega^{-*, *}(X, \End E)$.

\subsection{DGBV algebras from K\"{a}hler manifolds}
A K\"{a}hler manifold $X$ is automatically a symplectic manifold,
hence it one can obtain DGBV algebras from it as for a symplectic manifold.
Moreover,
as discovered in Cao-Zhou \cite{Cao-Zho1},
the type decomposition and the K\"{a}hler identities 
enable one to define two more DGBV algebras:
$(\Omega^{*, *}(X), \wedge, \bar{\partial}, -\sqrt{-1}\partial^*,
[\cdot\bullet\cdot]_{\-\sqrt{-1}\partial^*})$
and 
$(\Omega^{*, *}(X), \wedge, \partial, \sqrt{-1}\bar{\partial}^*,
[\cdot\bullet\cdot]_{\-\sqrt{-1}\bar{\partial}^*})$.

\section{$A_{\infty}$ algebraic mirror symmetry}

As in Kontsevich \cite{Kon},
we consider the mirror pair 
consists of a closed complex manifold $X$ 
and a closed symplectic manifold $\widehat{X}$.
The deformations of the complex structure on $X$ is described by
the Maurer-Cartan equation (see e.g. Kuranishi \cite{Kur})
$$\bar{\partial} \Gamma + \frac{1}{2}[\Gamma, \Gamma]_S = 0,$$
for $\Gamma \in \Omega^{-1, 1}(X)$.
We refer to the problem of solving this equation with 
$\Gamma \in \Omega^{-*, *}(X)$ as the extended deformation
problem of the complex structures.
This leads us to the differential G-algebra 
$(\Omega^{-*, *}(X), \wedge, \bar{\partial}, [\cdot, \cdot]_S)$.
Given a Hermitian metric,
Hodge theory provides a cohomological splitting of $\Omega^{-*, *}(X)$ (\S ?),
hence by \S \ref{sec:DiffGA},
one obtains an $A_{\infty}$ algebra structure 
and a canonically defined formal deformation of
this $A_{\infty}$ algebra structure on $H^{-*, *}(X)$.
The symplectic structure on $\widehat{X}$ determines a Poisson structure 
$w \in \Omega^{-2}(\widehat{X})$,
if $w + \gamma$ is another Poisson structure,
then from $[w + \gamma, w + \gamma]_S = 0$ one gets
$[w, \gamma]_S + \frac{1}{2}[\gamma, \gamma]_S = 0$.
Using the isomorphism $\Omega^{-*}(\widehat{X}) \to \Omega^*(\widehat{X})$,
one also gets a Maurer-Cartan equation
$$d \Gamma + \frac{1}{2}[\Gamma \bullet \Gamma]_{\Delta} = 0,$$
where $\Gamma \in \Omega^2(\widehat{X})$.
We refer to the problem of solving this equation with 
$\Gamma \in \Omega^*(X)$ as the extended deformation
problem of the Poisson structures.
This leads us to the DGBV algebra 
$(\Omega^*(\widehat{X}, \bC), \wedge, d, 
\Delta, [\cdot\bullet\cdot]_{\Delta})$.
Given a compatible Riemannian metric on $\widehat{X}$,
one then obtains
an $A_{\infty}$ algebra structure and 
a canonically defined formal deformation of
this $A_{\infty}$ algebra structure on $H^*(\widehat{X}, \bC)$.
The $A_{\infty}$ mirror symmetry means the existence of isomorphisms
of the above two $A_{\infty}$ algebras together 
with their canonically defined formal deformations.

If one considers topological open string theories,
then one encounters the deformation theories of Lagrangian submanifolds 
in symplectic manifolds and holomorphic vector bundles on complex manifolds
for $A$-model and $B$-model respectively.
(See e.g. Kontsevich \cite{Kon} and Vafa \cite{Vaf}.)
If $E$ is a holomorphic vector bundle on $X$,
the deformations of the holomorphic bundle structure on $E$ is
described by solving the following Maurer-Cartan equation:
$$\bar{\partial}_{\End E} \Gamma + \frac{1}{2}[\Gamma, \Gamma] = 0,$$
where $\Gamma \in \Omega^{0, 1}(X, \End E)$.
We refer to the problem of find solutions 
with $\Gamma$ in $\Omega^{0, *}(X, \End E)$ or $\Omega^{*, *}(X, \End E)$
or $\Omega^{-*, *}(X, \End E)$
as the extended deformation problem of the holomorphic bundle structure.
On these spaces we have differential Poisson algebra structures.
Using a Hermitian metric on $E$,
Hodge theory provides the cohomological splittings,
hence by \S \ref{sec:Poisson},
one obtains $A_{\infty}$ algebra structures 
and canonically defined deformations
on the the corresponding cohomology groups.
We conjecture that the similar constructions can be found for 
extended deformation problems of (special) Lagrangian submanifolds
so that one can formulate an $A_{\infty}$ algebraic mirror symmetry conjecture
as above.
It is interesting to examine other deformation problems 
in differential geometry and algebraic geometry in the same fashion.

We will compare the $A_{\infty}$ mirror symmetry for Calabi-Yau manifold
with the mirror symmetry formulated in terms of formal Frobenius manifolds
in Cao-Zhou \cite{Cao-Zho1, Cao-Zho4} in a separate paper.

\end{document}